# Comment

## Olivier Bousquet and Bernhard Schölkopf

Our contribution will be short, but we will try to compensate by being particularly opinionated. The field of support vector machines (SVMs) and related kernel methods has produced an impressive range of theoretical results, algorithms and success stories in real-world applications. While it originated in machine learning, it is also concerned with core problems of statistics and it is thus timely to publish a comprehensive article that discusses these methods from a statistician's point of view. We shall use this opportunity to make a few general comments, largely about the field rather than about the present paper.

Many papers about SVMs start off saying something like "SVMs are great because they are based on statistical learning theory" (this probably includes some of our own writings). Moguerza and Muñoz are more careful and only say that SVMs appeared in the context of statistical learning theory. What actually is the connection between SVMs and statistical learning theory?

Historically, SVMs and their precursors were (co-)developed by Vladimir Vapnik, one of the fathers of statistical learning theory. Statistical learning theory includes an analysis of machine learning which is independent of the distribution underlying the data. However, this analysis cannot provide any a priori guarantee that SVMs (or any other algorithm) will work well on a real-world problem. So what is special about SVMs, if anything?

In our view, what is special about SVMs is the combination of the following ingredients: first and foremost, the use of positive definite kernels; then regularization via the norm in the associated reproducing kernel Hilbert space; finally, the use of a convex loss function which is minimized by a classifier and not a regressor.

*Olivier Bousquet is Director of Research, Pertinence, F-75002 Paris, France e-mail: o.bousquet@pertinence.com. Bernhard Schölkopf is Professor and Director, Max Planck Institute for Biological Cybernetics, D-72076 Tübingen, Germany e-mail: bs@tuebingen.mpg.de.*

*The magic of kernels.* Positive definite kernels and their feature space interpretation do provide a very nice way to look at a whole class of algorithms; however, it is important to stress that they do not bring any *statistical* guarantee by themselves. The statistical guarantees available stem from the regularization (or learning theory) point of view. We shall return to this point below.

The main advantages of positive definite kernels are the following:

1. They allow easy construction of a nonlinear algorithm from a linear one, often without incurring additional computational cost.
2. They provide generality via the fact that they can be defined on nonvectorial data and do not, in general, require an explicit mapping to a reproducing kernel Hilbert space.

Historically, the first point was initially considered one of the major advantages of kernels and it triggered a significant number of kernel algorithms other than SVMs, starting with kernel principal component analysis (PCA). More recently, the second point has arguably taken over the role of the key selling point for kernel methods. The application of learning algorithms to nonvectorial data has become the field where nowadays a lot of the action is happening in the machine learning world, in particular concerning applications on structured data (e.g., in biology or natural language processing). We are curious to see whether the field of statistics will also embrace these possibilities.

*A sober look at the geometric interpretation.* The geometric point of view is an original way to look at SVMs and quite possibly the right way to come up with an algorithm like the SVM in the first place. However, it does not yield comprehensive statistical understanding. More precisely, there is no way to prove that large margin separating hyperplanes perform better than other types of hyperplanes independently of the distribution of the data.

Sure enough, the geometric point of view does provide intuition and motivates a large number of related algorithms, but one should not be fooled by geometric intuition or two-dimensional illustrations. The fact that data that are not linearly separable in input space suddenly becomes linearly separable in the so-called feature space (as depicted on Figure 1 of the main paper) has led to misconceptions. Indeed, the picture seems to suggest that the kernel has magically placed the two clouds of points in two separate regions of the space, and hence uncovered the right decision boundary.

The feature space often has a nonintuitive geometry. Let us take the example of the Gaussian Radial Basis Function (RBF) kernel. The corresponding feature space is of infinite dimension and the points are all mapped to the positive orthant of the unit sphere. *Any* two disjoint point sets in input space can be separated by a hyperplane in this feature space.

There is thus something mysterious happening in this space, but this space is but one way to look at things. We might instead directly look at the SVM algorithm and see that it loosely speaking tries to combine functions of the form $k(x_i, \cdot)$ using coefficients chosen to maximize the real-valued predictions $y_i f(x_i)$ on the training set. This brings us to the concept of *margin*. People usually say that maximizing the margin is good for generalization. There are two concepts of margin to be distinguished:

- *Geometric margin* (distance to the hyperplane). This is related to the norm of the weight vector, so that maximizing the margin corresponds to minimizing the norm (i.e., to regularization). Regularization can indeed lead to good generalization, provided the kind of smoothness enforced by the regularizer reflects the specifics of the problem.
- *Numerical margin* [i.e., the quantity $y_i f(x_i)$ which appears in the hinge loss used by the standard SVM]. The main reason why it makes sense to maximize this margin is because the hinge loss is a convex non-increasing upper bound of the classification loss, so that making $y_i f(x_i)$ large will ensure that the hinge loss is small and thus that we minimize the number of misclassification errors. However, this only means that minimizing the empirical hinge loss might lead to minimizing the empirical misclassification error, but does not guarantee that the expected misclassification error will be minimized as well.

These two notions are quite distinct, yet they are sometimes confused because they are entangled in the algorithms. For instance, if one minimizes the hinge loss over linear combinations of kernels and if there exists a combination such that the total hinge loss on the training set is zero, then this combination is not unique: we can multiply it by an arbitrary positive scale factor. Introducing a constraint on the norm of the weight vector is a natural way to remove this gauge freedom. This constraint is not innocent. It introduces a coupling between the numerical and the geometric margins: maximizing the geometric margin (in the context of an appropriate nonlinear kernel) leads to regularization which prevents overfitting by penalizing complex functions, while maximizing the numerical margin leads to minimization of the empirical error. Searching for a function with small empirical error while penalizing the complexity is the key to most reasonable learning algorithms.

*Convexity and loss functions.* Another attractive feature of positive definite kernels is that they allow nonlinearization of learning algorithm while preserving the *convexity* of the associated optimization problem. This is also one reason for the success of SVMs: the optimization problem is easier to handle than that of other algorithms such as artificial neural networks. The introduction of SVMs with kernels in the machine learning community suddenly moved the focus from optimization algorithms (e.g., multiple variants of gradient descent) to optimization criteria. This has created significant interest in convex functionals (for all kinds of problems such as model selection, semisupervised or unsupervised learning) and methods of convexifying existing functionals.

In the context of supervised learning, this search for convexity has led to the introduction of many different convex loss functions. However, something that has often been overlooked is the set of properties the loss function has to satisfy so that it leads to a consistent algorithm. For example, in the classification setting, a minimum requirement is that with sufficient data, minimizing the loss should lead to minimization of the misclassification error. For standard SVMs, the fact that the hinge loss satisfies this property was noticed relatively late (see reference [40] of the main paper) and, more surprisingly, in the context of multiclass classification this has been addressed only very recently. It has been proved in [1] that several variants of multiclass SVM do not have the required property. Of course, this is not to say that they perform poorly on a finite sample, but it is important to understand what an algorithm is aiming at and how it should behave as the sample size increases.

Moguerza and Muñoz are indeed aware of the fact that the minimizer of the hinge loss is the Bayes classifier (or rather is a function which has the same sign as that of the Bayes classifier), but they later say that there is still work to be done to provide a probabilistic interpretation of the output values produced by SVM classifiers. This is somewhat problematic because, at least asymptotically, there is no possible relationship between probabilities and output values. [This follows from the consistency property: With an appropriate kernel, the values of the function produced by the SVM algorithm will converge to exactly $+1$ or $-1$ on all the points where $P(Y|X) \in ]0, 1/2[ \cup ]1/2, 1[$ so that the value of $f(X)$ will have no relationship to $P(Y|X)$.] Hence, on a finite sample, if a relationship occurs, it will likely be by pure chance (or because the kernel happens to regularize exactly in the way needed for the preferred functions to look like the conditional probability density function).

To conclude this section with a more philosophical viewpoint, let us mention that the SVM algorithm also reinforces the belief that one should be concerned about the objective rather than about the model: what is important is not whether one can identify the "true" target function; rather, one should try to find *some* function, from a large class, which will perform well. This belief is shared by many researchers in the machine learning community, and it probably distinguishes them from "classical" statisticians, as argued, for example, in [3].

*Theoretical considerations.* Regarding the statistical analysis of the SVM algorithm, besides the works cited in the paper, there are a few additional references that are worth mentioning; for example, [6] first proved universal consistency of $L_1$-SVM with a Gaussian kernel, while Steinwart and Scovel [8] and Steinwart [7] obtained rates of convergence under various conditions. Also, more recently, the consistency of SVM has been proved by Vert and Vert [9] in the case where the regularization parameter is held fixed, but the kernel width goes to zero. This suggests that there is a coupling between both types of regularization (provided by a small norm of the function and a large kernel width).

It is now clear that the VC dimension is not the right parameter to capture the rates of convergence, especially when studying real-valued functions classes. Alternative possibilities (based on Rademacher averages) along with finite-sample performance bounds can be found, for example, in [2].

Progress has also been made in understanding the role of sparsity in SVMs. First of all, the number of support vectors is asymptotically linear in the sample size if the Bayes error is nonzero. Second, on large data bases the number of support vectors is usually too large for fast testing (hence the development of reduced set methods which can be applied to nonsparse models [4, 5]).

*Why do SVM work so well in practice?* There is probably no theoretical answer to this question. The fact that they are universally consistent is surely interesting, but does not explain anything about finite sample performance on real-world data sets (e.g., the $k$-nearest neighbor algorithm is also universally consistent). The sparsity also does not explain it. Regularization (by the kernel width and by the function norm) surely plays a role (by preventing overfitting) but this cannot be quantified. Indeed, in statistical terms, one can only tell the effect of regularization on the *variance* but not on the *bias*, at least if one does not make specific assumptions on the smoothness of the target function. The only possible answer to this question might thus be that on those problems where SVMs excel, the kernel that is used induces a regularizer that incorporates appropriate prior knowledge about the problems or, equivalently, it captures the right notion of similarity. In a large majority of applications, the Gaussian RBF kernel is used and its success simply means that the Euclidean distance in input space is locally meaningful for those problems. [Indeed, the Gaussian kernel incorporates a notion of similarity which is a monotonic function of the Euclidean distance. In this case, the SVM produces a "local rule": The prediction at a given point is a weighted combination of the labels of nearby points (where the weight mainly depends on the distance and is adapted by the coefficients $\lambda_i$ which appear in equation (3.4) of the main paper).]

*Future directions for research.* Although, as explained by Moguerza and Muñoz, the SVM algorithm in itself has several interesting merits, we think that what is most important about it is its impact on the field of machine learning and statistics. It has introduced new concepts and ideas that have considerably influenced their progress, and we expect that the acquired momentum will lead to further advances, in domains such as structured learning, joint kernels (mixing inputs and outputs), links to graphical models, and semisupervised learning, to name but a few. In a different direction, one could try to extend the notion of kernel so as to handle higher level similarities, such as analogies (which can be considered as similarities between pairs of examples).

There are also several important questions that need to be addressed so as to bridge the gap between basic research and applications. For instance, there is no satisfactory method for choosing the parameters other than using cross-validation, which can be an obstacle in applications. Moreover, there are still significant computational issues arising from the implementation of SVM-like algorithms using nonlinear kernels for large-scale problems.